\documentclass[12pt]{amsart}

\usepackage[english]{babel}
\usepackage[T1]{fontenc}
\usepackage[utf8]{inputenc}
\usepackage{lmodern} 

\usepackage{amsmath}
\usepackage{amstext}
\usepackage{amsfonts}
\usepackage{amssymb}
\usepackage{amsthm}
\usepackage{mathtools}
\usepackage[colorinlistoftodos]{todonotes}

\usepackage{verbatim}
\usepackage{float}

\usepackage[colorlinks=true, allcolors=blue]{hyperref}
\usepackage{theoremref} 

\newtheorem{theorem}{Theorem}[section]
\newtheorem{lemma}[theorem]{Lemma}
\newtheorem{prop}[theorem]{Proposition}
\newtheorem{cor}[theorem]{Corollary}

\newtheorem*{cor*}{Corollary}
\newtheorem*{thm*}{Theorem}
\newtheorem*{lem*}{Lemma}

\newtheorem*{prop*}{Proposition}

\theoremstyle{definition}
\newtheorem{definition}[theorem]{Definition}

\newtheorem{example}[theorem]{Example}
\newtheorem*{defn*}{Definition}

\newtheorem{remark}[theorem]{Remark}

\newcommand{\Ad}{\operatorname{Ad}}

\newcommand{\bC}{{\mathbb{C}}}
\newcommand{\bE}{{\mathbb{E}}}

\newcommand{\Sub}{\operatorname{Sub}}

\newcommand{\C}{\mathbb{C}} 

\DeclareRobustCommand{\SkipTocEntry}[5]{}

\newcommand{\R}{\mathbb{R}}  
\newcommand{\Z}{\mathbb{Z}}
\newcommand{\N}{\mathbb{N}}

\newcommand{\norm}[1]{ \left\| #1 \right\| }

\newcommand{\setbuilder}[2] { \left\{ #1 \enskip \middle| \enskip #2 \right\} }

\newcommand{\abs}[1]{ \left| #1 \right| }

\newcommand{\closure}[1]{\overline{#1}}

\newcommand{\set}[1]{\left\{#1\right\}}

\newcommand{\boundary}{\partial}
\newcommand{\defeq}{\coloneqq}

\newcommand{\wtilde}[1]{\widetilde{#1}}

\newcommand{\injectsinto}{\hookrightarrow}
\newcommand{\surjectsonto}{\twoheadrightarrow}

\renewcommand{\Re}{\operatorname{Re}}

\DeclareMathOperator{\weakstar}{w*}

\title[Powers averaging for commutative crossed products]{A generalized Powers averaging property for commutative crossed products}

\author[T. Amrutam]{Tattwamasi Amrutam}
\address{Department of Mathematics\\ University of Houston\\USA}
\email{tamrutam@math.uh.edu}
\author[D. Ursu]{Dan Ursu}
\address{Department of Pure Mathematics, University of Waterloo, 200 University Avenue
West, Waterloo, Ontario, N2L 3G1, Canada}
\email{dursu@uwaterloo.ca}

\thanks{Second author supported by the Natural Sciences and Engineering Research Council of Canada (NSERC) [grant number PGSD3-535032-2019]. Deuxi\`{e}me auteur financ\'{e} par le Conseil de recherches en sciences naturelles et en g\'{e}nie du Canada (CRSNG) [num\'{e}ro de subvention PGSD3-535032-2019].}
\begin{document}

\begin{abstract}
    We prove a generalized version of Powers' averaging property that characterizes simplicity of reduced crossed products $C(X) \rtimes_\lambda G$, where $G$ is a countable discrete group, and $X$ is a compact Hausdorff space which $G$ acts on minimally by homeomorphisms. As a consequence, we generalize results of Hartman and Kalantar on unique stationarity to the state space of $C(X) \rtimes_\lambda G$ and to Kawabe's generalized space of amenable subgroups $\Sub_a(X,G)$. This further lets us generalize a result of the first named author and Kalantar on simplicity of intermediate C*-algebras. We prove that if $C(Y) \subseteq C(X)$ is an inclusion of unital commutative $G$-C*-algebras with $X$ minimal and $C(Y) \rtimes_\lambda G$ simple, then any intermediate C*-algebra $A$ satisfying $C(Y) \rtimes_\lambda G \subseteq A \subseteq C(X) \rtimes_\lambda G$ is simple.
\end{abstract}

\maketitle

\tableofcontents

\section{Introduction}
The notion of Powers' averaging property for discrete groups has played an important role in recent years in questions about simplicity related to reduced group C*-algebras and reduced crossed products. In this paper, we introduce a generalized version of Powers' averaging property for reduced crossed products of the form $C(X) \rtimes_\lambda G$, and prove that it is equivalent to simplicity of the crossed product. We then derive various consequences.

First, we recall the notion of Powers' averaging property, along with a brief history of recent applications. Let $G$ be a countable discrete group, and let $P(G)$ denote the set of probability measures on $G$. For convenience, we will denote the finitely supported probability measures on $G$ by $P_f(G)$. Recall that we canonically have an action of $P(G)$ on any $G$-C*-algebra $A$ as follows: given $\mu \in P(G)$ and $a \in A$, we let
\[ \mu a = \sum_{g \in G} \mu(g) (g \cdot a). \]

The group $G$ is said to be C*-simple if its reduced group C*-algebra $C^*_\lambda(G)$ is simple. It was shown independently in \cite[Theorem~6.3]{Kennedy} and \cite[Theorem~4.5]{Haagerup} that C*-simplicity is equivalent to an averaging property originally considered by Powers, which can most conveniently be stated as follows: $G$ is said to have \textit{Powers' averaging property} if for any $a \in C^*_\lambda(G)$, we have
\[ \tau_\lambda(a) \in \closure{\setbuilder{\mu a}{\mu \in P_f(G)}}. \]
Here, $\tau_\lambda$ denotes the canonical trace on $C^*_\lambda(G)$, and we are canonically viewing $\C \subseteq C^*_\lambda(G)$. The set $P_f(G)$ above can be replaced by $P(G)$ instead. It is also clear that it suffices to check only the $a \in C^*_\lambda(G)$ satisfying $\tau_\lambda(a) = 0$, as it is always possible to ``normalize'' an arbitrary $a \in C^*_\lambda(G)$ by considering $a - \tau_\lambda(a)$.

It was later shown by Hartman and Kalantar in the proof of \cite[Theorem~5.1]{HartKal} that averaging with respect to all of $P_f(G)$ is not necessary, and that if $G$ is countable, then Powers' averaging property for $C^*_\lambda(G)$ is equivalent to the existence of a single measure $\mu \in P(G)$ (not guaranteed to have finite support) satisfying $\mu^n a \to \tau_\lambda(a)$ for all $a \in C^*_\lambda(G)$.

A generalization of Powers' averaging property to reduced (twisted) crossed products of unital C*-algebras and C*-simple groups was given by 
\cite[Section~3]{KB}. They showed that the same averaging property holds for elements $a \in A \rtimes_\lambda G$ satisfying $\mathbb{E}(a) = 0$, where $\mathbb{E} : A \rtimes_\lambda G \to A$ denotes the canonical conditional expectation, i.e.
\[ 0 \in \closure{\setbuilder{\mu a}{\mu \in P_f(G)}}. \]

The ideas mentioned above were used by the first named author and Kalantar in \cite{AK} to show that if $G$ is C*-simple and the action of $G$ on a compact Hausdorff space $X$ is minimal, then not only is the reduced crossed product $C(X) \rtimes_\lambda G$ simple, but so is any intermediate C*-algebra lying between $C^*_\lambda(G)$ and $C(X) \rtimes_\lambda G$.

Of course, if $G$ is not C*-simple, then $A \rtimes_\lambda G$ can never have Powers' averaging property, but it is still possible for the crossed product to be simple. An easy example is $C(\mathbb{T}) \rtimes_\lambda \Z$, where $\Z$ acts on the circle $\mathbb{T}$ by an irrational rotation. For this reason, we introduce a generalized version of Powers' averaging which does turn out to be equivalent to simplicity in the end. Let $X$ be a compact Hausdorff space equipped with an action of $G$ by homeomorphisms.

We define in Section~\ref{section_generalized_probability_measures}, and in particular in \thref{definition_generalized_probability_measure}, the spaces $P(G,C(X))$ and $P_f(G,C(X))$ of what we call \textit{generalized $(G,C(X))$-probability measures}. Given an inclusion of unital $G$-C*-algebras $C(X) \subseteq A$, the space $P(G,C(X))$ canonically admits a left action on $A$, and a right action on the state space $S(A)$. With this, we are able to conveniently generalize Powers' averaging property to crossed products as follows:

\begin{definition}
\thlabel{definition_generalized_powers_averaging}
    Let $G$ be a countable discrete group acting on a compact Hausdorff space $X$ by homeomorphisms, and let $\bE : C(X) \rtimes_\lambda G \to C(X)$ denote the canonical conditional expectation. We say that $C(X) \rtimes_\lambda G$ has the \textit{generalized Powers' averaging property} if for every $a \in C(X) \rtimes_\lambda G$ with $\mathbb{E}(a) = 0$, we have
    \[ 0 \in \closure{\setbuilder{\mu a}{\mu \in P_f(G,C(X))}}. \]
\end{definition}

One might define other generalized analogues of Powers' averaging property, for example requiring that $\mathbb{E}(a)$ lie in the above set given any $a \in C(X) \rtimes_\lambda G$ not necessarily satisfying $\mathbb{E}(a) = 0$. It is not immediately obvious, however, that this is equivalent with the version in \thref{definition_generalized_powers_averaging}, as considering $a - \mathbb{E}(a)$ for an arbitrary $a \in C(X) \rtimes_\lambda G$ just tells us that for any $\varepsilon > 0$, there is a $\mu \in P_f(G,C(X))$ with the property that $\norm{\mu a - \mu \mathbb{E}(a)} < \varepsilon$. However, unlike in the case where $C(X) = \C$, we do not have that $\mu \mathbb{E}(a) = \mathbb{E}(a)$ in general. Hence, our first main result is perhaps a bit surprising:

\begin{theorem}
\thlabel{simple_iff_powers_averaging}
    Let $G$ be a countable discrete group acting on a compact Hausdorff space $X$ by homeomorphisms, and assume that the action is minimal. Let $\bE : C(X) \rtimes_\lambda G \to C(X)$ denote the canonical conditional expectation. The following are equivalent:
    \begin{enumerate}
        \item $C(X) \rtimes_\lambda G$ is simple.
        \item $C(X) \rtimes_\lambda G$ has the generalized Powers' averaging property.
        \item $\bE(a) \in \closure{\setbuilder{\mu a}{\mu \in P_f(G,C(X))}}$ for all $a \in C(X) \rtimes_\lambda G$.
        \item $\nu(\bE(a)) \in \closure{\setbuilder{\mu a}{\mu \in P_f(G,C(X))}}$ for all $a \in C(X) \rtimes_\lambda G$ and $\nu \in P(X)$.
    \end{enumerate}
\end{theorem}

Next we generalize Hartman and Kalantar's results. It is worth noting that they operate under a slightly different action of $P(G)$ on any $G$-C*-algebra, with a convolution product given by
\[ \mu * a = \sum_{g \in G} \mu(g) (g^{-1} \cdot a) \]
for any $\mu \in P(G)$, and a left action of $P(G)$ on $S(A)$ given by
\[ \mu * \phi = \sum_{g \in G} \mu(g) (g \cdot \phi). \]
However, this is only a minor technicality to keep in mind, and it is easy to rephrase their results (which we do) in terms of the actions we use in our paper.

As previously mentioned, they show that Powers' averaging property for $C^*_\lambda(G)$ is equivalent to the existence of a measure $\mu \in P(G)$ with the property that $\mu^n a \to \tau_\lambda(a)$ for any $a \in C^*_\lambda(G)$ \cite[Theorem~5.1]{HartKal}. As a consequence, the only state $\phi \in S(C^*_\lambda(G))$ that is \textit{$\mu$-stationary} (that is, satisfying $\phi \mu = \phi$) is the canonical trace $\tau_\lambda$, and this is in fact a characterization of C*-simplicity of $G$ \cite[Theorem~5.2]{HartKal}. Similar result holds in the crossed product setting:

\begin{theorem}
\thlabel{onemeasureconv}
    Let $G$ be a countable discrete group acting on a compact Hausdorff space $X$ by homeomorphisms, and assume that the action is minimal. Let $\bE : C(X) \rtimes_\lambda G \to C(X)$ denote the canonical conditional expectation. If $C(X) \rtimes_\lambda G$ is simple, then there is a generalized measure $\mu \in P(G,C(X))$ with the property that $\mu^n a \to 0$ whenever $\bE(a) = 0$. Optionally, we may also require that $\mu$ have full support.
\end{theorem}

\begin{cor}
\thlabel{statstatesandsimplicity}
    Let $G$ be a countable discrete group acting on a compact Hausdorff space $X$ by homeomorphisms, and assume that the action is minimal. Let $\bE : C(X) \rtimes_\lambda G \to C(X)$ denote the canonical conditional expectation. Then the crossed product $C(X) \rtimes_\lambda G$ is simple if and only if there is some $\mu \in P(G,C(X))$ with full support and with the property that any $\mu$-stationary state on $C(X) \rtimes_\lambda G$ is of the form $\nu \circ \bE$ for some $\mu$-stationary $\nu \in P(X)$.
\end{cor}

It is worth noting in \thref{statstatesandsimplicity} that, given a generalized measure $\mu \in P(G,C(X))$, there is no guarantee that there be a \textit{unique} $\mu$-stationary measure $\nu \in P(X)$. If one could strengthen the averaging in \thref{simple_iff_powers_averaging} (4) to work with a single measure, for example if there were some $\nu \in P(X)$ and $\mu \in P(G,C(X))$ such that $\mu^n a \to \nu(\bE(a))$ for any $a \in C(X) \rtimes_\lambda G$, then it would be possible to obtain uniqueness of $\nu$ as well. However, we were unable to prove such a result.

Our first application of Powers' averaging property is a natural generalization of the main result in \cite{AK}.

\begin{theorem}
\thlabel{simplcityofintmalgebras}
    Let $G$ be a countable discrete group, and assume that $C(Y) \subseteq C(X)$ is an equivariant inclusion of commutative unital $G$-C*-algebras. Assume moreover that the action of $G$ on $X$ is minimal. If $C(Y)\rtimes_{\lambda}G$ is simple, then every intermediate C*-algebra $A$ satisfying $C(Y) \rtimes_\lambda G \subseteq A \subseteq C(X) \rtimes_\lambda G$ is simple.
\end{theorem}

One other result of Hartman and Kalantar that we generalize is the following. Denote the space of amenable subgroups of $G$ by $\Sub_a(G)$. This is naturally a compact Hausdorff space if we equip it with the Chabauty topology, which is the topology induced by viewing this canonically as a subset of $2^G$ (the power set of $G$), and it also carries a $G$-action by homeomorphisms given by conjugation. It has been known for a few years that the dynamics on this space characterizes C*-simplicity, with \cite[Theorem~4.1]{Kennedy} essentially stating that $G$ is C*-simple if and only if $\set{\set{e}}$ is the unique minimal component in $\Sub_a(G)$, and \cite[Corollary~5.7]{HartKal} stating that C*-simplicity is equivalent to unique stationarity of $\delta_{\set{e}}$ with respect to some $\mu \in P(G)$.

The dynamical analogue for crossed products $C(X) \rtimes_\lambda G$ (where $X$ is minimal) is a result of Kawabe \cite[Theorem~6.1]{Kawabe}. Consider the space
\[ \Sub_a(X,G) \coloneqq \setbuilder{(x,H)}{x \in X, H \leq G_x, \text{ and } H \text{ amenable}},\]
where $G_x$ denotes the stabilizer subgroup of $x$. This is again a compact Hausdorff space with $G$-action given by $s \cdot (x,H) = (sx, sHs^{-1})$, and Kawabe's result amounts to saying that $C(X) \rtimes_\lambda G$ is simple if and only if the only minimal component in $\Sub_a(X,G)$ is $X \times \set{\set{e}}$. This hints that there should also be a ``unique stationarity result'' involving measures supported on this minimal component.

\begin{cor}
\thlabel{unique_stationarity_generalized_amenable_subgroup_space}
    Let $G$ be a countable discrete group acting on a compact Hausdorff space $X$ by homeomorphisms, and assume that the action is minimal. Let $\Sub_a(X,G)$ denote Kawabe's generalized space of amenable subgroups, and view $C(X) \subseteq C(\Sub_a(X,G))$ as dual to the canonical projection $\Sub_a(X,G) \surjectsonto X$ mapping $(x,H)$ to $x$. The crossed product $C(X) \rtimes_\lambda G$ is simple if and only if there is some $\mu \in P(G,C(X))$ with the property that any $\mu$-stationary measure in $P(\Sub_a(X,G))$ is supported on $X \times \set{\set{e}}$.
\end{cor}

\addtocontents{toc}{\SkipTocEntry}
\section*{Acknowledgements}
The authors owe a huge debt of gratitude to Mehrdad Kalantar and Matthew Kennedy, who are their advisors respectively, for many useful discussions surrounding the problem.
The first named author is also grateful to Narutaka Ozawa for many helpful discussions about the averaging property. The authors also thank Sven Raum, Eusebio Gardella, Shirley Geffen and Yongle Jiang for taking the time to carefully read a near complete draft of this paper and giving helpful feedback. The authors would also like to thank the anonymous referee for the detailed reading of our paper and for their comments and suggestions which enhanced the exposition of the paper.

\addtocontents{toc}{\SkipTocEntry}
\subsection*{Recent Development} Upon completion of this paper, Narutaka Ozawa sent us his preprint where he has shown (via different techniques) that, in the non-minimal case, a slight variation on the generalized Powers' averaging property is equivalent to the action of $G$ on $X$ having the residual intersection property. We remark that for minimal $G$-spaces $X$, residual intersection property is equivalent to the crossed product $C(X) \rtimes_\lambda G$ being simple.

\section{The space of generalized probability measures}
\label{section_generalized_probability_measures}

To establish notation, $G$ will denote a countable discrete group, and $X$ will denote a compact Hausdorff space which $G$ acts on by homeomorphisms. All C*-algebras and morphisms are assumed to be unital.

We define the notion of a generalized probability measure. As motivation, consider the case of a $G$-C*-algebra $A$. Given a fixed $a \in A$, any probability measure $\mu \in P(G)$ provides a convenient way of representing a convex combination of the elements $\setbuilder{g \cdot a}{g \in G}$. Namely, we may define $\mu a \defeq \sum_{g \in G} \mu(g) g \cdot a$.

With this in mind, we want a space of generalized probability measures which represents \textit{C(X)-convex combinations}. For convenience, we first review this notion here:

\begin{definition}
\thlabel{definition_CX_convex}
    Assume $C(X) \subseteq A$ is an inclusion of unital C*-algebras, and let $K \subseteq A$. We say that $K$ is \textit{$C(X)$-convex} if, given finitely many $f_1, \dots, f_n \in C(X)$ with $\sum_{i=1}^n f_i^2 = 1$, and any $a_1, \dots, a_n \in K$, we have $\sum_{i=1}^n f_i a_i f_i \in K$. Such a sum is called a \textit{$C(X)$-convex combination} of $a_1, \dots, a_n$.
\end{definition}

\begin{remark}
    The usual notion of $C(X)$-convex combinations is slightly more general, and deals with sums of the form $\sum_{i=1}^n f_i^* a_i f_i$, where $\sum_{i=1}^n f_i^*f_i = 1$ and $f_i$ is no longer assumed to be positive. For our purposes, we will stick with the definition given in \thref{definition_CX_convex}, as working with positive $f_i$ is in particular necessary for \thref{ineq-av-1} later on.
\end{remark}

If $C(X) \subseteq A$ is an inclusion of unital $G$-C*-algebras, and $a \in A$, we want our notion of generalized probability measures to represent $C(X)$-convex combinations of $\setbuilder{g \cdot a}{g \in G}$.

\begin{definition}
\thlabel{definition_generalized_probability_measure}
    Consider a formal sum
    \[ \mu = \sum_{s \in G} \sum_{i \in I_s} f_i s f_i \]
    (where all $I_s$ are disjoint), with the properties $f_i \geq 0$, $f_i \neq 0$, and $\sum_{s \in G} \sum_{i \in I_s} f_i^2 = 1$.
    Equivalenty, we may also combine the above double-sum into a single sum
    \[ \mu = \sum_{i \in I} f_i s_i f_i \]
    if we allow repetition among the group elements $s_i$. We say that $\mu$ is a \textit{generalized $(G,C(X))$-probability measure}, and denote the set of all such generalized measures by $P(G,C(X))$. The set of all finite-sum generalized measures is denoted by $P_f(G,C(X))$. Given a unital $G$-C*-algebra $A$ containing an equivariant copy of $C(X)$, and any $\mu = \sum_{i \in I} f_i s_i f_i \in P(G,C(X))$, we have a unital and completely positive map on $A$ given by
	\[ \mu a = \sum_{i \in I} f_i (s_i \cdot a) f_i. \]
	Moreover, this induces a right action on the state space $S(A)$, given by $(\phi \mu)(a) = \phi(\mu a)$.
\end{definition}

\begin{remark}
		\thlabel{remark_repetition_of_group_elements}
		We note a couple of things. First, observe that for a fixed $s \in G$, it is in general not possible to simplify an expression of the form $f_1 s f_1 + f_2 s f_2$ as a single $hsh$ for some $h \in C(X)$. This is a consequence of noncommutativity. If $C(X) \subseteq A$ and we were to consider the action of this element on some $a \in A$, this becomes
		$$ f_1 (s \cdot a) f_1 + f_2 (s \cdot a) f_2, $$
		which in general is not equal to any $h (s \cdot a) h$ for any $h \in C(X)$. Because of this, repetition of group elements is allowed, and this is also reasoning behind the choice of terminology and notation, namely ``finite-sum'' and $P_f(G,C(X))$, as opposed to ``compactly supported'' and $P_c(G,C(X))$. Given a generalized probability measure $\mu = \sum_{i \in I} f_i s_i f_i$, it is possible to have infinitely many elements $i \in I$ with $s_i$ all being equal. In other words, it is possible to have an infinite sum (which doesn't necessarily simplify to a finite sum) that is still ``compactly supported'' on $G$.
\end{remark}
	
The rest of this section is dedicated to proving various technicalities and basic properties of the space $P(G,C(X))$. First, the following is an easy exercise in functional analysis:

\begin{lemma}
\thlabel{lemma_banach_space_infinite_sum_convergence}
    Assume $X$ is a Banach space, and $\sum_{i \in I} x_i$ is an infinite unordered sum. Then this sum converges if and only if for all $\varepsilon > 0$, there exists a finite set $F \subseteq I$ such that for all finite sets $J \subseteq I \setminus F$, we have $\norm{\sum_{j \in J} x_j} < \varepsilon$.
\end{lemma}

From this, we obtain two important results:

\begin{cor}
    Any $\mu = \sum_{i \in I} f_i s_i f_i$ in $P(G,C(X))$ has the property that $I$ is countable. In particular, $P(G,C(X))$ is indeed a set.
\end{cor}

\begin{proof}
    Consider the sum $\sum_{i \in I} f_i^2 = 1$ and $\varepsilon = \frac{1}{n}$ ($n \in \N$) in \thref{lemma_banach_space_infinite_sum_convergence}. Then necessarily, only finitely many $f_i^2$ can have norm at least $\frac{1}{n}$. Hence, at most countably many $f_i$ can be nonzero.
\end{proof}

\begin{cor}
    Assume $C(X) \subseteq A$ is an inclusion of $G$-C*-algebras, $\mu = \sum_{i \in I} f_i s_i f_i \in P(G,C(X))$, and $a \in A$. Then the sum given by $\sum_{i \in I} f_i (s_i \cdot a) f_i$ is convergent, or in other words, the value $\mu a$ is well-defined. Moreover, the map $a \mapsto \mu a$ is unital and completely positive.
\end{cor}

\begin{proof}
    We first prove this for positive $a$. Let $\varepsilon > 0$, and let $F \subseteq I$ be such that for all finite $J \subseteq I \setminus F$, we have $\norm{\sum_{j \in J} f_j^2} < \varepsilon$. Then
    \[ \norm{\sum_{j \in J} f_j (s_j \cdot a) f_j} \leq \norm{a} \norm{\sum_{j \in J} f_j^2} < \varepsilon \norm{a}. \]
    To see that general values of $a$ work, let $\mu_F = \sum_{i \in F} f_i s_i f_i$ for finite $F \subseteq I$. Writing $a$ as a finite linear combination of four positive elements $\sum_{k=1}^4 c_k a_k$, we have that each of the nets $(\mu_F a_k)_F$ converges for each $k$. In particular, the net
    \[ \mu_F a = \sum_{k=1}^4 c_k \mu_F a_k \]
    must therefore also be convergent. The fact that $a \mapsto \mu a$ is completely positive follows from the fact that $a \mapsto \mu_F a$ is completely positive for each finite $F \subseteq I$.
\end{proof}

Similar to how $P(G)$ is a convex semigroup, we have that the spaces $P_f(G,C(X))$ and $P(G,C(X))$ also form semigroups, and moreover satisfy an appropriate notion of $C(X)$-convexity.

\begin{prop}
\thlabel{generalized_probability_measures_are_CX_convex}
    The space $P_f(G,C(X))$ is $C(X)$-convex, in the sense that given finitely many $\set{g_j}_{j \in J} \subseteq C(X)$ with $\sum_{j \in J} g_j^2 = 1$ and any $\set{\mu_j}_{j \in J} \subseteq P_f(G,C(X))$ with $\mu_j = \sum_{i \in I_j} f_i s_i f_i$ (and all index sets $I_j$ disjoint for distinct values of $j$), we have that
    \[ \sum_{j \in J} g_j \mu_j g_j \defeq \sum_{j \in J} \sum_{i \in I_j} g_j f_i s_i f_i g_j \]
    also lies in $P_f(G,C(X))$. The same is true for $P(G,C(X))$, except that $J$ can be infinite.
\end{prop}

\begin{proof}
    We prove the case of $P(G,C(X))$, as the case of $P_f(G,C(X))$ is almost the same except without needing to worry about convergence. Observe that, given any finite subsets $F \subseteq J$ and $F_j \subseteq I_j$ for $j \in F$, we have
    \[ \sum_{j \in F} \sum_{i \in F_j} g_j^2 f_i^2 = \sum_{j \in F} g_j^2 \sum_{i \in F_j} f_i^2 \leq \sum_{j \in F} g_j^2 \cdot 1 \leq 1. \]
    Moreover, given $\varepsilon > 0$, if we choose $F \subseteq J$ finite with $\sum_{j \in J} g_j^2 \geq 1 - \varepsilon$ and finite $F_j \subseteq I_j$ for $j \in F$ with $\sum_{i \in F_j} f_i^2 \geq 1 - \varepsilon$, then
    \[ \sum_{j \in F} \sum_{i \in F_j} g_j^2 f_i^2 = \sum_{j \in F} g_j^2 \sum_{i \in F_j} f_i^2 \geq \sum_{j \in F} g_j^2 \cdot (1 - \varepsilon) \geq (1 - \varepsilon)^2. \]
    This proves that $\sum_{j \in J} \sum_{i \in I_j} g_j^2 f_i^2 = 1$.
\end{proof}

It is perhaps worthwhile to do an example of a $C(X)$-convex combination.

\begin{example}
    Consider $X = [0,1]$, let $G$ be some group acting on $X$, let $s$ and $t$ be distinct group elements, and consider the following set of generalized probability measures and coefficients:
    \begin{align*}
		\mu_1 &= \sqrt{x} s \sqrt{x} + \sqrt{1-x} t \sqrt{1-x} \\
		\mu_2 &= \sqrt{1-x} s \sqrt{1-x} + \sqrt{x} t \sqrt{x} \\
		g_1 &= \frac{1}{\sqrt{2}} \\
		g_2 &= \frac{1}{\sqrt{2}}
	\end{align*}
	We have that
	\begin{align*}
		g_1 \mu_1 g_1 + g_2 \mu_2 g_2 &= (\frac{1}{\sqrt{2}} \sqrt{x}) s (\frac{1}{\sqrt{2}} \sqrt{x}) + (\frac{1}{\sqrt{2}} \sqrt{1-x}) s (\frac{1}{\sqrt{2}} \sqrt{1-x}) \\
		&+ (\frac{1}{\sqrt{2}} \sqrt{x}) t (\frac{1}{\sqrt{2}} \sqrt{x}) + (\frac{1}{\sqrt{2}} \sqrt{1-x}) t (\frac{1}{\sqrt{2}} \sqrt{1-x})
	\end{align*}
	which cannot be simplified further (see \thref{remark_repetition_of_group_elements}), and is left as-is, with group elements being repeated in the expression.
\end{example}

\begin{remark}
\thlabel{remark_mu_a_smallest_G_invariant_CX_convex_subset}
    Given an inclusion of $G$-C*-algebras $C(X) \subseteq A$, $a \in A$, finitely many $\set{g_j}_{j \in J} \subseteq C(X)$ with $\sum_{j \in J} g_j^2 = 1$, and $\set{\mu_j}_{j \in J} \subseteq P_f(G,C(X))$, we have that
    \[ \left(\sum_{j \in J} g_j \mu_j g_j\right)(a) = \sum_{j \in J} g_j \mu_j(a) g_j. \]
    Consequently,
    \[ \setbuilder{\mu a}{\mu \in P_f(G,C(X))} \]
    is $C(X)$-convex as well. In fact, it is the smallest $G$-invariant, $C(X)$-convex subset of $A$ containing $a$.
\end{remark}

Now we define a semigroup structure on $P_f(G,C(X))$ and $P(G,C(X))$.

\begin{prop}
\thlabel{generalized_probability_measures_are_semigroup}
    The space $P(G,C(X))$ is a semigroup under the following multiplication: given $\mu = \sum_{i \in I} f_i s_i f_i$ and $\nu = \sum_{j \in J} g_j t_j g_j$, let
    \[ \mu \nu \defeq \sum_{i \in I} \sum_{j \in J} (f_i (s_i g_j)) (s_i t_j) ((s_i g_j) f_i). \]
    Moreover, $P_f(G,C(X))$ is a subsemigroup of $P(G,C(X))$.
\end{prop}

\begin{proof}
    Observe that, given any finite subsets $F_I \subseteq I$ and $F_J \subseteq J$, we have
    \[ \sum_{i \in F_I} \sum_{j \in F_J} (f_i^2 (s_i g_j))^2 = \sum_{i \in F_I} f_i^2 s_i \left(\sum_{j \in F_J} g_j^2\right) \leq \sum_{i \in F_I} f_i^2 s_i 1 \leq 1. \]
    Moreover, any finite subset of $I \times J$ is contained in a set of the form $F_I \times F_J$. Finally, given $\varepsilon > 0$, if one chooses $F_I$ and $F_J$ to be such that $\sum_{i \in F_I} f_i^2 \geq 1 - \varepsilon$ and $\sum_{j \in F_J} g_j^2 \geq 1 - \varepsilon$, then we have
    \[ \sum_{i \in F_I} \sum_{j \in F_J} (f_i (s_i g_j))^2 = \sum_{i \in F_I} f_i s_i \left(\sum_{j \in F_J} g_j^2\right) \geq \sum_{i \in F_I} f_i s_i(1-\varepsilon) \geq (1-\varepsilon)^2. \]
    This proves that $\sum_{i \in I} \sum_{j \in J} (f_i (s_i g_j))^2 = 1$, and so this multiplication on $P(G,C(X))$ is well-defined. Associativity is tedious but not hard to check. The fact that $P_f(G,C(X))$ is a subsemigroup is clear.
\end{proof}

\begin{remark}
\thlabel{remark_generalized_probability_measures_semigroup_action}
    The multiplication on $P(G,C(X))$ is defined in such a way so that if $C(X) \subseteq A$ is an inclusion of unital $G$-C*-algebras, $\mu_1,\mu_2 \in P(G,C(X))$, and $a \in A$, then
    \[ (\mu_1 \mu_2)(a) = \mu_1(\mu_2 a). \]
    In other words, we canonically have a left semigroup action of $P(G,C(X))$ on $A$, and consequently a right semigroup action on $S(A)$.
\end{remark}

Let $C(X) \subseteq A$ be an inclusion of unital $G$-C*-algebras, and let $\mu \in P(G,C(X))$ be a generalized measure. We say that a state $\phi \in S(A)$ is \textit{$\mu$-stationary} if $\phi \mu = \phi$, and denote the set of all $\mu$-stationary states on $A$ by $S_{\mu}(A)$. Observe that this definition makes sense even for C*-subalgebras that don't necessarily contain $C(X)$, but are at least $\mu$-invariant. It is not hard to see that $\mu$-stationary states always exist and can be extended to a larger C*-algebra. The proof is a mere modification of \cite[Proposition~4.2]{HartKal}. We include it for the sake of completeness.

\begin{prop}
\thlabel{existext}
Suppose that $C(X) \subseteq A$ is an inclusion of unital $G$-C*-algebras, $\mu \in P(G,C(X))$, and $B \subseteq A$ is a $\mu$-invariant unital C*-subalgebra. Then every $\mu$-stationary state $\tau \in S_\mu(B)$ can be extended to a $\mu$-stationary state $\eta \in S_\mu(A)$. In particular, $S_\mu(A)$ is always nonempty.
\begin{proof}
Let $K = \setbuilder{\zeta \in S(A)}{\zeta|_B = \tau}$, a compact convex set. For any $\mu \in P(G,C(X))$, the map $\Phi_{\mu}:K\to K$ defined by $\Phi_{\mu}(\zeta)=\zeta\mu$ is an affine continuous map. It is well-known that $\Phi_{\mu}$ has a fixed point, say $\eta$. Then, $\eta \in S_{\mu}(A)$ and $\eta|_B=\tau$. To see that $S_\mu(A)$ is nonempty, let $B = \C$.
\end{proof}
\end{prop}

We wish to define an appropriate notion of full support for measures in $P(G,C(X))$. For this, the following observation will come in useful.

\begin{lemma}
    Assume $X$ is a Banach space, and $\sum_{i \in I} x_i$ is an infinite unordered sum that converges in norm. Then given any $J \subseteq I$, the sum $\sum_{j \in J} x_j$ also converges.
\end{lemma}

\begin{proof}
    We know by \thref{lemma_banach_space_infinite_sum_convergence} that given any $\varepsilon > 0$, there is a finite subset $F \subseteq I$ such that for any finite subset $E \subseteq I \setminus F$, we have $\norm{\sum_{i \in E} x_i} < \varepsilon$. But then, letting $F' = F \cap J$, it is clear that for any finite set $E' \subseteq J \setminus F'$, we also have $\norm{\sum_{j \in E'} x_j} < \varepsilon$.
\end{proof}

\begin{definition}
\thlabel{fullsupport}
We say that a generalized measure $\mu\in P(G,C(X))$ has full support if, writing
\[\mu=\sum_{s\in G} \sum_{i \in I_s} f_i s f_i,\]
we have that for each $s \in G$, there is some $\delta > 0$ such that $\sum_{i \in I_s} f_i^2 \geq \delta$. Equivalently (by compactness of $X$), given any $s \in G$ and $x \in X$, we can find $i \in I_s$ such that $f_i(x)>0$.
\end{definition}

\section{Proof of generalized Powers averaging}
In this section, we prove \thref{simple_iff_powers_averaging}. To give a brief overview, we first recall how this is proven in the case of the usual reduced group C*-algebra.

Both in \cite{Kennedy} and \cite{Haagerup}, which independently prove C*-simplicity is equivalent to Powers' averaging property, the key tool used was the dynamical characterization of C*-simplicity of $G$ in terms of its action on the Furstenberg boundary $\partial_F G$ (see \cite[Theorem~1.5]{kalantar_kennedy_boundaries} or \cite[Theorem~1.1]{breuillard_kalantar_kennedy_ozawa_c_simplicity}). The Furstenberg boundary of a group was originally developed by Furstenberg \cite{Furstenberg1973} (see also \cite{Furstenberg}) as a topological object, but it can also be realized as the $G$-injective envelope of the complex numbers $\bC$ (see \cite[Definition~3.1, Theorem~3.11]{kalantar_kennedy_boundaries}). Briefly, we recall the topological characterization below:

\begin{definition}
    Let $X$ be a $G$-space. A measure $\nu \in P(X)$ is called \textit{contractible} if $\setbuilder{\delta_x}{x \in X} \subseteq \closure{G \nu}^{\weakstar}$. A \textit{$G$-boundary} is a $G$-space $X$ with the additional property that every measure $\nu \in P(X)$ is contractible.
\end{definition}

It is worth noting that, from the perspective of convexity, being a $G$-boundary is the same as saying that $P(X)$ is irreducible as a compact convex $G$-space.

\begin{theorem}
    There is a universal $G$-boundary $\boundary_F G$ (known as the Furstenberg boundary of $G$) with the property that every $G$-boundary $X$ is the image of $\boundary_F G$ under a surjective, $G$-equivariant map.
\end{theorem}

It is shown in \cite{Kennedy} that C*-simplicity of $G$ is equivalent to $\set{\tau_\lambda}$ being the unique $G$-boundary in $S(C^*_\lambda(G))$. A similar result can consequently be achieved in the entire dual space of $C^*_\lambda(G)$, and a Hahn-Banach separation argument yields that this is equivalent to Powers' averaging property.

There is a generalized notion of boundaries introduced in \cite[Section~7]{KS19}, which is used for dealing with noncommutative crossed products $A \rtimes_\lambda G$. However, this notion is more technical, as it involves working with matrix convex sets and matrix state spaces. It is possible to develop a similar notion, but using only usual convex sets instead of matrix convex ones, and use this in the case of commutative crossed products $C(X) \rtimes_\lambda G$. However, \cite{Kawabe}, which deals with proving equivalences of simplicity of such crossed products, does not develop such a theory of generalized boundaries. Instead, the generalized notion of boundary necessary here is developed in \cite[Section~3]{Naghavi}, albeit from the perspective of compact sets and their measures rather than from convex sets.

From Naghavi's results, we use (part of) \cite[Theorem~3.2]{Naghavi} and the discussion following it, which we quickly paraphrase as follows:

\begin{theorem}
\thlabel{naghavithm}
    For a countable discrete group $G$, let $X$ be a minimal $G$-space, and let $C(X) \subseteq C(Y)$ be an inclusion of unital $G$-C*-algebras. The following are equivalent:
    \begin{enumerate}
        \item $C(Y)$ is a $G$-essential extension of $C(X)$.
        \item Given any measure $\nu \in P(Y)$ with the property that the restriction $\nu|_{C(X)} \in P(X)$ is contractible, we have that $\nu$ is also contractible.
    \end{enumerate}
    In particular, the above is true for the $G$-injective envelope $C(Y) = I_G(C(X))$, which is a maximal $G$-essential extension of $C(X)$.
\end{theorem}

\begin{definition}
    Assume $G$ is a countable discrete group, $X$ is a minimal $G$-space, and $C(X) \subseteq C(Y)$ is an inclusion of unital $G$-C*-algebras. We say that $Y$ is a \textit{$(G,X)$-boundary} if $C(Y)$ satisfies any of the equivalent conditions above. Furthermore, we let $\boundary_F(G,X)$ denote the Gelfand spectrum of the $G$-injective envelope of $C(X)$, i.e. $I_G(C(X)) = C(\boundary_F(G,X))$. It can be shown that $\boundary_F(G,X)$ is universal among all $(G,X)$-boundaries.
\end{definition}

We also recall the result \cite[Theorem~3.4]{Kawabe} of Kawabe characterizing the ideal intersection property for $C(X) \rtimes_\lambda G$ in the special case of minimal dynamical systems (in which case, simplicity and the intersection property coincide).

\begin{theorem}
\thlabel{kawabethm}
    For a countable discrete group $G$, let $X$ be a minimal $G$-space, and let $\partial_F(G,X)$ be the Gelfand spectrum of the $G$-injective envelope of $C(X)$. Then the following are equivalent:
    \begin{enumerate}
        \item The crossed product $C(\partial_F(G,X)) \rtimes_\lambda G$ is simple.
        \item The crossed product $C(X) \rtimes_\lambda G$ is simple.
        \item The action of $G$ on $\partial_F(G,X)$ is free.
    \end{enumerate}
\end{theorem}

With the above, we are ready to tackle proving our generalized version of Powers' averaging property. The following lemma shows that, although measures on minimal spaces need not be contractible in general, they have the weaker property that arbitrary measures can still be pushed to Dirac masses using $P(G,C(X))$.

\begin{lemma}
\thlabel{lemma_generalized_contractibility}
Assume $X$ is a minimal $G$-space, and fix any $x \in X$. There is a net $(\mu_\lambda) \subseteq P_f(G,C(X))$ with the property that for any $\nu \in P(X)$, we have $\nu \mu_\lambda \xrightarrow{\weakstar} \delta_x$.
\end{lemma}

\begin{proof}
    Fix an open neighbourhood $V$ of $x$. Observe that, given any $y \in X$, we have that $Gy$ is dense in $X$ by minimality. In particular, there exists some $s \in G$ with the property that $sy \in V$, or equivalently, $y \in s^{-1} V$. It follows that the sets $sV$ form an open cover of $X$, and so there is some finite subcover $s_1 V, \dots, s_n V$. Now let $F_1, \dots, F_n \in C(X)$ be a partition of unity subordinate to this open cover, let $f_i = F_i^{1/2}$, and let
    \[ \mu_V = \sum_{i=1}^n f_i s_i f_i. \]
    It is not hard to see that, given any $\nu \in P(X)$, $\nu \mu_V$ is a measure with support contained in $\closure{V}$. It follows that the net $(\nu \mu_V)$, indexed by open neighbourhoods of $x$ ordered under reverse inclusion, converges weak* to $\delta_x$.
\end{proof}

This allows us to push arbitrary measures towards the trivial boundary in $S(C(X) \rtimes_\lambda G)$ in the case of simple crossed products.

\begin{prop}
\thlabel{generalized_averaging_state_space}
    Let $X$ be a minimal $G$-space, and assume that the crossed product $C(X) \rtimes_\lambda G$ is simple. Then given any state $\phi \in S(C(X) \rtimes_\lambda G)$, we have that
    \[ \setbuilder{\nu \circ \bE}{\nu \in P(X)} \subseteq \closure{\setbuilder{\phi \mu}{\mu \in P_f(G,C(X))}}^{\weakstar}. \]
\end{prop}

\begin{proof}
    For convenience, denote the latter set above by $K$. By $G$-invariance, convexity, and weak*-closure, it suffices to prove that $K$ contains $\delta_x \circ \bE$ for any single point $x \in X$.
    
    To this end, let $I_G(C(X)) = C(\boundary_F(G,X))$ be the $G$-injective envelope of $C(X)$. Extend the state $\phi$ to a state $\wtilde{\phi}$ on $C(\boundary_F(G,X)) \rtimes_\lambda G$. By \thref{lemma_generalized_contractibility}, we can find a net $(\mu_\lambda) \subseteq P_f(G,C(X))$ with the property that $\wtilde{\phi}|_{C(X)} \mu_\lambda \to \delta_x$ for some $x \in X$. Dropping to a subnet if necessary, we have that $\wtilde{\phi} \mu_\lambda \to \psi \in S(C(\boundary_F(G,X)) \rtimes_\lambda G)$ with the property that $\psi|_{C(X)} = \delta_x$. Observe that $\psi|_{C(X) \rtimes_\lambda G} \in K$.
    
    Minimality tells us that $\psi|_{C(X)}$ is contractible, and so $\psi|_{C(\boundary_F(G,X))}$ is contractible as well by \thref{naghavithm}. This tells us that there is a net $(g_i)$ with $\psi|_{C(\boundary_F(G,X))} g_i \to \delta_y$ for some $y \in \boundary_F(G,X)$. Again dropping to a subnet yields a state $\eta \in S(C(\boundary_F(G,X)) \rtimes_\lambda G)$ with the property that $\eta|_{C(\boundary_F(G,X))} = \delta_y$. Observe once more that $\eta|_{C(X) \rtimes_\lambda G} \in K$.
    
    We claim that $\eta|_{C(X) \rtimes_\lambda G}$ is the state we are looking for. Simplicity of $C(X) \rtimes_\lambda G$ implies that the action of $G$ on $\boundary_F(G,X)$ is free by \thref{kawabethm}. From here, the rest is a common argument. We know that $C(\boundary_F(G,X))$ lies in the multiplicative domain of $\eta$. Thus, it suffices to show that $\eta(\lambda_t) = 0$ for $t \neq e$, as this will imply that for any $f \in C(\boundary_F(G,X))$, we have
    \[ \eta(f \lambda_t) = f(y) \eta(\lambda_t) = 0. \]
    Now let $f \in C(\boundary_F(G,X))$ be such that $f(y) = 1$ and $f(ty) = 0$. This is possible because $ty \neq y$. We have
    \[ f(y) \eta(\lambda_t) = \eta(f \lambda_t) = \eta(\lambda_t (t^{-1} \cdot f)) = \eta(\lambda_t) f(ty). \]
    This forces $\eta(\lambda_t) = 0$, as desired.
\end{proof}

The fact that arbitrary functionals on a C*-algebra are a finite linear combination of states gives us a similar result on the entire dual space $(C(X) \rtimes_\lambda G)^*$.

\begin{prop}
\thlabel{generalized_averaging_dual_space}
    Let $X$ be a minimal $G$-space, and assume that the crossed product is simple. Then given any $\omega \in (C(X) \rtimes_\lambda G)^*$, we have
    \[ \setbuilder{\omega(1) \nu \circ \bE}{\nu \in P(X)} \subseteq \closure{\setbuilder{\omega \mu}{\mu \in P_f(G,C(X))}}^{\weakstar}. \]
\end{prop}

\begin{proof}
    Again, for convenience, denote this latter set by $K$. Write $\omega = \sum_{i=1}^4 c_i \phi_i$, a linear combination of four states. By \thref{generalized_averaging_state_space}, we can find a net $(\mu_\lambda) \subseteq P_f(G,C(X))$ with the property that $\phi_1 \mu_\lambda \to \nu_1 \circ \bE$. Dropping to a subnet if necessary, we may also assume that $(\phi_i \mu_\lambda)$ are all convergent to some $\phi_i'$ for $i \geq 2$. In particular, $(\omega \mu_\lambda)$ converges to some $\omega' \in K$ with the property that $\omega' = c_1 \nu_1 \circ \bE + \sum_{i=2}^4 c_i \phi_i'$. Noting that the set $\setbuilder{\nu \circ \bE}{\nu \in P(X)}$ is weak*-closed and closed under the right action of $P_f(G,C(X))$, repeating this averaging trick three more times nets us (no pun intended) that there is some element in $K$ of the form $\eta \circ \bE$ satisfying $\eta(1) = \omega(1)$.
    
    Writing $\eta = \sum_{i=1}^4 d_i \psi_i$ as a linear combination of four states on $C(X)$, fixing $x \in X$, and letting $(\mu_j)$ be as in \thref{lemma_generalized_contractibility}, we have that $(\eta \circ \bE) \mu_j = (\eta \mu_j) \circ \bE$ converges to $\omega(1) \delta_x \circ \bE$, which must lie in $K$. Minimality of $X$ and $G$-invariance, weak*-closure, and convexity of $K$ yield that every $\omega(1) \nu \circ \bE$ lies in $K$ as well.
\end{proof}

From here, it is an application of the Hahn-Banach separation argument that gives us the strong generalized Powers' averaging property. Conversely, lack of nontrivial ideals can be directly deduced even from just being able to average elements $a \in C(X) \rtimes_\lambda G$ satisfying $\bE(a) = 0$.

\begin{proof}[Proof of \thref{simple_iff_powers_averaging}]
    $(1) \implies (4)$ Given that the extreme points of $P(X)$ are the Dirac masses $\delta_x$, it suffices to prove the following: if $a \in C(X) \rtimes_\lambda G$ and $x \in X$, then
    \[ \bE(a)(x) \in \closure{\setbuilder{\mu a}{\mu \in P_f(G,C(X))}}. \]
    Assume otherwise, so that there is some $a \in C(X) \rtimes_\lambda G$ and $x \in X$ for which this doesn't hold. Then there is some functional $\omega \in (C(X) \rtimes_\lambda G)^*$ and $\alpha \in \R$ with the property that
    \[ \Re \omega(\bE(a)(x)) < \alpha \leq \Re \omega(\mu a) \quad \forall \mu \in P_f(G,C(X)). \]
    However, given that $\omega(\bE(a)(x)) = \omega(1) \bE(a)(x)$, and by \thref{generalized_averaging_dual_space}, $\omega(\mu a)$ can be made arbitrarily close to $\omega(1) (\delta_x \circ \bE)(a)$, this cannot happen.
    
    $(4) \implies (3)$ Our aim is to show that we may approximate $\bE(a)$ by $C(X)$-convex combinations of $\bE(a)(x)$, where $x \in X$. Let $a \in C(X) \rtimes_\lambda G$, and let $\varepsilon > 0$. Given any $x \in X$, by continuity of $\bE(a) \in C(X)$, there is some open neighbourhood $U_x$ of $x$ for which $\abs{\bE(a)(x) - \bE(a)(y)} < \varepsilon$ for all $y \in U_x$. By compactness, there is some finite subcover $U_{x_1}, \dots, U_{x_n}$ of $X$. Let $F_i$ be a partition of unity subordinate to the open cover, and let $f_i = F_i^{1/2}$. Observe that, given any $x \in X$, we have
    \begin{align*}
        &\abs{\sum_{i=1}^n f_i(x) (\bE(a)(x_i)) f_i(x) - \bE(a)(x)} \\
        &=\abs{\sum_{i=1}^n f_i(x) (\bE(a)(x_i) - \bE(a)(x)) f_i(x)} \\
        &\leq \sum_{i=1}^n f_i(x)^2 \abs{\bE(a)(x_i) - \bE(a)(x)} \\
        &< \sum_{i=1}^n f_i(x)^2 \cdot \varepsilon \\
        &= \varepsilon
    \end{align*}
    Thus, we have $\norm{\sum_{i=1}^n f_i \bE(a)(x_i) f_i - \bE(a)} < \varepsilon$. By \thref{remark_mu_a_smallest_G_invariant_CX_convex_subset}, we have that $\bE(a) \in \closure{\setbuilder{\mu a}{\mu \in P_f(G,C(X))}}$, as this set is closed under $C(X)$-convex combinations.
    
    $(3) \implies (2)$ This direction is clear.
    
    $(2) \implies (1)$ Let $I$ be any nontrivial ideal of $C(X) \rtimes_\lambda G$, and let $a$ be any nonzero element of $I$. Replacing $a$ by $a^*a$, we may assume without loss of generality that $a$ is a nonzero positive element. Faithfulness of the canonical expectation tells us that $\mathbb{E}(a)$ is also a nonzero positive element, and so there is some $\varepsilon > 0$ and open subset $U \subseteq X$ with the property that $\bE(a)(x) > \varepsilon$ for all $x \in U$. Using the same trick as in the proof of \thref{lemma_generalized_contractibility}, minimality of $X$ gives us that $X = s_1 U \cup \dots \cup s_n U$ for finitely many $s_i \in G$. Hence, replacing $a$ by $s_1 a + \dots + s_n a$, we may assume without loss of generality that $\bE(a) > \varepsilon$. If we choose $\mu \in P(G,C(X))$ so that $\norm{\mu(a - \bE(a))} < \frac{\varepsilon}{2}$, then as this value is in particular self-adjoint, we have that $\mu(a - \bE(a)) \geq - \frac{\varepsilon}{2}$. Consequently,
    \[ \mu a = \mu(\bE(a)) + \mu (a - \bE(a)) \geq \varepsilon - \frac{\varepsilon}{2} = \frac{\varepsilon}{2}. \]
    In particular, $\mu a \in I$ is invertible, which gives us that $I$ is the entire crossed product $C(X) \rtimes_\lambda G$.
\end{proof}

\section{Unique stationarity and applications}

This section generalizes the various results in \cite{HartKal} on equivalence between C*-simplicity and unique stationarity of the canonical trace in $C^*_\lambda(G)$, along with its consequences.

It is worth noting that one cannot expect simplicity of $C(X) \rtimes_\lambda G$ to be equivalent to unique stationarity of an element of $S(C(X) \rtimes_\lambda G)$, even with respect to a generalized measure $\mu \in P(G,C(X))$. This is because of the fact that there may not exist a uniquely stationary state on $C(X)$, and any $\mu$-stationary state on $C(X)$ will extend to one on the whole crossed product. The natural fix is to instead expect that the $\mu$-stationary states on $C(X) \rtimes_\lambda G$ all be of the form $\nu \circ \bE$, where $\nu$ ranges over the $\mu$-stationary measures $\nu \in P(X)$. It is also worth noting that one cannot expect to work with the usual notion of measure $\mu \in P(G)$, as this would again imply unique stationarity of $\tau_\lambda \in S(C^*_\lambda(G))$. However, this is equivalent to C*-simplicity of $G$ \cite[Theorem~5.2]{HartKal}, which is by no means necessary for the crossed product $C(X) \rtimes_\lambda G$ to be simple - take, for example, $C(\mathbb{T}) \rtimes_\lambda \Z$, where $\Z$ acts on the circle $\mathbb{T}$ by an irrational rotation.\\

We begin with the observation that averaging elements in the reduced group C*-algebra $C^*_\lambda(G)$, even with respect to a generalized measure $\mu \in P(G,C(X))$, is enough to average elements in the crossed product $C(X) \rtimes_\lambda G$.

\begin{lemma}
\thlabel{ineq-av-1}
Let $X$ be a minimal $G$-space, and let $\mu \in P(G,C(X))$. Then given any $t \in G$ and $f \in C(X)$, we have
\[ \norm{\mu (f \lambda_t)} \leq \norm{f} \norm{\mu \lambda_t}. \]
\end{lemma}

\begin{proof}
It is well-known that the crossed product $\ell^\infty(G) \rtimes_\lambda G$ (the uniform Roe algebra), can canonically be viewed as a C*-subalgebra of $B(\ell^2(G))$. Fixing $x_0 \in X$ gives us a unital $G$-equivariant injective *-homomorphism $\iota : C(X) \injectsinto \ell^\infty(G)$, given by $\iota(f)(t) = f(tx_0)$. This lets us view $C(X) \rtimes_\lambda G$ as a C*-subalgebra of $B(\ell^2(G))$ as well.

Write $\mu = \sum_{i \in I} g_i s_i g_i$. Given $\xi \in \ell^2(G)$ and $r \in G$, we have
\begin{align*}
    ((\mu (f \lambda_t))\xi)(r)
    &= \left( \sum_{i \in I} g_i (s_i f) (s_i t s_i^{-1} g_i) \lambda_{s_i t s_i^{-1}} \xi \right)(r) \\
    &= \sum_{i \in I} g_i(rx_0) f(s_i^{-1} r x_0) g_i(s_i t^{-1} s_i^{-1} r x_0) \xi(s_i t^{-1} s_i^{-1} r).
\end{align*}
Now letting $\abs{\xi} \in \ell^2(G)$ be given by $\abs{\xi}(r) = \abs{\xi(r)}$, we note that $\norm{\abs{\xi}} = \norm{\xi}$. Moreover, we have
\begin{align*}
    \norm{(\mu (f \lambda_t))\xi}^2 &= \sum_{r \in G} \abs{\sum_{i \in I} g_i(rx_0) f(s_i^{-1} r x_0) g_i(s_i t^{-1} s_i^{-1} r x_0) \xi(s_i t^{-1} s_i^{-1} r)}^2 \\
    &\leq \norm{f}^2 \sum_{r \in G} \left( \sum_{i \in I} g_i(rx_0) g_i(s_i t^{-1} s_i^{-1} r x_0) \abs{\xi(s_i t^{-1} s_i^{-1} r)} \right)^2 \\
    &= \norm{f}^2 \norm{(\mu \lambda_t) \abs{\xi}}^2.
\end{align*}
It follows that $\norm{\mu(f \lambda_t)} \leq \norm{f} \norm{\mu \lambda_t}$.
\end{proof}

It is also an easy remark that Powers' averaging property can be made to work with finitely many elements at once.

\begin{lemma}
\thlabel{lemma_powers_averaging_multiple_elements}
    Assume $C(X) \rtimes_\lambda G$ has Powers' averaging property. Then given any $a_1, \dots, a_n \in C(X) \rtimes_\lambda G$ satisfying $\bE(a_i) = 0$, and $\varepsilon > 0$, there is some $\mu \in P(G,C(X))$ with the property that $\norm{\mu a_i} < \varepsilon$ for all $i = 1, \dots, n$.
\end{lemma}

\begin{proof}
    Let $\mu_1 \in P(G,C(X))$ be such that $\norm{\mu_1 a_1} < \varepsilon$. Choosing $\mu_{k+1} \in P(G,C(X))$ inductively by letting $\mu_{k+1}$ be such that $\norm{\mu_{k+1} (\mu_k \dots \mu_1 a_{k+1})} < \varepsilon$, we see that $\mu = \mu_n \dots \mu_1$ is the generalized measure we are looking for.
\end{proof}

\begin{proof}[Proof of \thref{onemeasureconv}]
    First, we claim that there is such a measure that works for all elements $a \in C^*_\lambda(G) \subseteq C(X) \rtimes_\lambda G$ satisfying $\tau_\lambda(a) = 0$. This is a near-verbatim repeat of the proof of \cite[Theorem~5.1]{HartKal}.
    We repeat the construction of $\mu$ here, along with the appropriate modifications.
    
    Let $(n_k)$ be an increasing sequence of positive integers satisfying $\left(\sum_{i=1}^k \frac{1}{2^i}\right)^{n_k} < \frac{1}{2^k}$, and let $(a_i)$ be any dense sequence in the unit ball of $\ker \tau_\lambda \subseteq C^*_\lambda(G)$. Let $\mu_1 \in P(G,C(X))$ be anything. Using \thref{lemma_powers_averaging_multiple_elements}, we may inductively build $\mu_l$ for $l \geq 2$ so that
    \[ \norm{\mu_l \mu_{k_r} \dots \mu_{k_1} a_s} < \frac{1}{2^l} \]
    for all $1 \leq s,k_1,\dots,k_r < l$ and $0 \leq r < n_l$. Here, by $r = 0$, we mean that $\mu_l \mu_{k_r} \dots \mu_{k_1} a_s$ becomes $\mu_l a_s$. A tedious computation then shows that $\mu = \sum_{l=1}^\infty \frac{1}{2^l} \mu_l$ will satisfy $\mu^n a \to 0$ for any $a \in \ker \tau_\lambda$.
    
    To force $\mu$ to have full support, observe that if $(s_n)_{n \in \N}$ is an enumeration of $G$, then the measure
    \[ \nu = \sum_{n=1}^\infty \frac{1}{2^{n+1}} s_n \frac{1}{2^{n+1}} \in P(G,C(X)) \]
    has full support. Then fixing any $l$ and letting $\alpha > 0$ be sufficiently small, we may replace $\mu_l$ by $\alpha \nu + (1-\alpha) \mu_l$ and still satisfy the required approximation properties above. Thus, without loss of generality, some $\mu_l$ has full support, and hence so does $\mu$.
    
    Now, to see that $\mu^n a \to 0$ whenever $a \in C(X) \rtimes_\lambda G$ satisfies $\bE(a) = 0$, we first prove this for elements $a_0 = f_1 \lambda_{t_1} + \dots + f_n \lambda_{t_n}$, where $t_i \neq e$. Note that by \thref{ineq-av-1}, we have
    \[ \norm{\mu^n a_0} \leq \sum_{i=1}^n \norm{\mu (f_i \lambda_{t_i})} \leq \sum_{i=1}^n \norm{f_i} \norm{\mu \lambda_{t_i}} \to 0. \]
    Now given an arbitrary $a$ with $\bE(a) = 0$, and $\varepsilon > 0$, we may choose $a_0$ as before with $\norm{a - a_0} < \varepsilon$. Choosing $N$ such that, given $n \geq N$, we have $\norm{\mu^n a_0} < \varepsilon$, we also have
    \[ \norm{\mu^n a} \leq \norm{\mu^n a_0} + \norm{\mu^n (a-a_0)} < \varepsilon + \varepsilon = 2 \varepsilon. \]
\end{proof}

\begin{remark}
    In the above proof, if we instead wanted to directly construct a generalized measure $\mu \in P(G,C(X))$ with the property that $\mu^n a \to 0$ for all $a \in C(X) \rtimes_\lambda G$ with $\bE(a) = 0$, as opposed to $a \in C^*_\lambda(G)$ with $\tau_\lambda(a) = 0$, we would have required separability of $\ker \bE \subseteq C(X) \rtimes_\lambda G$, which requires separability of $C(X)$ (metrizability of $X$). Proceeding with $\ker \tau_\lambda \subseteq C^*_\lambda(G)$ first and then lifting the averaging to the entire crossed product avoids this extra assumption. There are natural examples of spaces on which $G$ acts that are not metrizable. For example, if $G$ is not amenable, then the Furstenberg boundary $\boundary_F G$ is such a space \cite[Corollary~3.17]{kalantar_kennedy_boundaries}.
\end{remark}

For a minimal $G$-space $X$, it is well known that if $\mu \in P(G)$ has full support, then any $\mu$-stationary state on $C(X)$ is faithful. A similar result holds for generalized probability measures (with the definition of full support given in \thref{fullsupport}).

\begin{lemma}
\thlabel{faithfullandfullsupport}
Let $X$ be a minimal $G$-space. Let $\mu \in P(G,C(X))$ be a generalized probability measure with full support. Then every $\mu$-stationary state on $C(X)$ is faithful.
\begin{proof}
Let $f \in C(X)$ be such that $f \geq 0$ and $f \neq 0$. Let $\mu = \sum_{s \in G} \sum_{i \in I_s} f_i s f_i$ be a generalized probability measure with full support. Since $f$ is nonzero, there exists $x_0\in X$ such that $f(x_0)>0$. It follows from $X$ being minimal that, for every $x\in X$, there exists $s_x \in G$ such that $f(s_x^{-1}x)>0$. Moreover, since $\mu$ has full support, there exists $i_x \in I_{s_x}$ such that $f_i(x)>0$.
Therefore,
\[ \mu(f)(x) = \sum_{s\in G} \sum_{i \in I_s} f_i(x) f(s^{-1}x) f_i(x) > f_{i_x}(x) f(s_x^{-1}x) f_{i_x}(x) > 0. \]
By compactness of $X$, it follows that there exists a $\delta>0$ such that $\mu(f) \geq \delta$.
Consequently, for any $\mu$-stationary state $\tau$ on $C(X)$, we see that
\[ \tau(f) = \tau(\mu(f)) \geq \delta. \]
Hence, $\tau$ is faithful.
\end{proof}

\end{lemma}
\begin{proof}[Proof of \thref{statstatesandsimplicity}] Let $X$ be a minimal $G$-space.
Suppose that $C(X)\rtimes_{\lambda}G$ is simple. Let $\mu \in P(G,C(X))$ be the generalized measure obtained from \thref{onemeasureconv} and let $\tau$ be a $\mu$-stationary state on $C(X)\rtimes_{\lambda}G$. Then, for any $a\in C(X)\rtimes_{\lambda}G$ with $\bE(a) = 0$, we have that
\[ \tau(a) = \tau(\mu^n a) \to \tau(0) = 0, \]
and so for general $a \in C(X) \rtimes_\lambda G$, we have
\[ \tau(a) = \tau(\bE(a)) + \tau(a - \bE(a)) = \tau(\bE(a)). \]
In other words, $\tau = \tau|_{C(X)} \circ \bE$.

On the other hand, suppose that there exists a generalized probability measure $\mu \in P(G,C(X))$ with full support along with the property that every $\mu$-stationary state on $C(X) \rtimes_\lambda G$ is of the form $\nu \circ \bE$ for some $\mu$-stationary $\nu \in P(X)$. By faithfulness of $\bE$ and \thref{faithfullandfullsupport}, every $\mu$-stationary state on $C(X)\rtimes_{\lambda}G$ is faithful. This is enough to guarantee that the $C(X) \rtimes_\lambda G$ is simple - the proof is similar to \cite[Proposition~4.9]{HartKal}.

Assume that there was a nontrivial ideal $I \subseteq C(X) \rtimes_\lambda G$. Observe that the quotient map $\pi : C(X) \rtimes_\lambda G \to (C(X) \rtimes_\lambda G)/I$ is nonfaithful, as nontrivial ideals always contain nonzero positive elements. Moreover, the quotient $(C(X) \rtimes_\lambda G)/I$ is canonically a $G$-C*-algebra, with $t \in G$ acting by $\Ad \pi(\lambda_t)$, and the quotient map $\pi$ is $G$-equivariant. In particular, we still canonically have $C(X) \subseteq (C(X) \rtimes_\lambda G)/I$ (under the quotient map $\pi$) by minimality of $X$, and so there is at least one $\mu$-stationary state $\phi \in S_\mu((C(X) \rtimes_\lambda G)/I)$ by \thref{existext}. The composition $\phi \circ \pi$ is a $\mu$-stationary state on $C(X) \rtimes_\lambda G$ that is not faithful, contradicting our earlier conclusion.
\end{proof}

One should notice that $G$-simplicity doesn't necessarily pass to sub-algebras and therefore, simplicity for invariant sub-algebras of simple crossed products shouldn't be expected to hold in general. Consider, for example, any simple C*-algebra $A$, any $C*$-simple group $G$ acting on $A$ trivially, and any abelian $C*$-subalgebra $B \subseteq A$. However, given an inclusion of unital $G$-C*-algebras $C(Y) \subset C(X)$ (via a factor map $\pi:X \to Y$), since any $G$-invariant C*-subalgebra $A$, $C(Y) \subset A \subset C(X)$ is of the form $C(Z)$ where $Z$ is an equivariant factor of $X$, and minimality passes to factors, it follows from the characterization of Kawabe \cite[Theorem~6.1]{Kawabe} that $C(Z) \rtimes_\lambda \Gamma$ is simple. We follow arguments similar to the proof of \cite[Theorem~1.3]{AK} to deal with general intermediate C*-subalgebras between $C(Y) \rtimes_\lambda G$ and $C(X) \rtimes_\lambda G$, not necessarily of the above form.

\begin{proof}[Proof of \thref{simplcityofintmalgebras}]
    By \thref{onemeasureconv}, there exists a generalized measure $\mu \in P(G,C(Y))$ with full support and the property that $\mu^n a \to 0$ whenever $a \in C(Y) \rtimes_\lambda G$ is such that $\bE(a) = 0$. Observe that we canonically have $P(G,C(Y)) \subseteq P(G,C(X))$. Since $X$ is minimal and $\mu$ has full support, it follows from \thref{faithfullandfullsupport} that every $\mu$-stationary state $\nu$ on $C(X)$ is faithful, and since $\bE$ is also faithful, it follows that every $\mu$-stationary state on $C(X) \rtimes_\lambda G$ of the form $\nu \circ \bE$ is faithful. We claim that the proof is complete once we establish that every $\mu$-stationary state $\tau$ on $C(X)\rtimes_{\lambda}G$ is of the form $\nu \circ \mathbb{E}$.
    
    Indeed, if this is the case, let $A$ be any intermediate C*-algebra of the form $C(Y) \rtimes_\lambda G \subseteq A \subseteq C(X) \rtimes_\lambda G$. Suppose that $I$ is a proper closed two-sided ideal of $A$. Then the action of $G$ on $A$ induces an action of $G$ on $A/I$ (as $I$ is necessarily $G$-invariant). Moreover, by minimality of $Y$, we also canonically have $C(Y) \subseteq A/I$. By \thref{existext}, there exists a $\mu$-stationary state $\varphi$ on $A/I$. Upon composing $\varphi$ with the canonical quotient map $A \to A/I$, we obtain a $\mu$-stationary state $\wtilde{\varphi}$ on $A$ which vanishes on $I$. Using \thref{existext} again, extend $\wtilde{\varphi}$ to a $\mu$-stationary state $\tau$ on $C(X) \rtimes_{\lambda} G$. By our assumption, $\tau$ being of the form $\nu \circ \bE$, is faithful. But $\tau|_I = 0$, which cannot occur if $I$ is nontrivial.
    
    We return to the question of showing that every $\mu$-stationary state on $C(X) \rtimes_\lambda G$ is indeed of the form $\nu \circ \bE$ for some $\mu$-stationary $\nu \in P(X)$. We claim that $\mu^n a \to 0$ whenever $a \in C(X) \rtimes_\lambda G$ (not just $C(Y) \rtimes_\lambda G$) satisfies $\bE(a) = 0$. To see this, first let $f \in C(X)$ and $t \neq e$. \thref{ineq-av-1} tells us that
    \[ \norm{\mu^n(f \lambda_t)} \leq \norm{f} \norm{\mu^n \lambda_t} \to 0. \]
    It follows that for finite linear combinations $a_0 = f_1 \lambda_{t_1} + \dots + f_n \lambda_{t_n}$, where $f_i \in C(X)$ and $t_i \neq e$, we have $\mu^n a_0 \to 0$ as well. Finally, let $a \in C(X) \rtimes_\lambda G$ with $\bE(a) = 0$, $\varepsilon > 0$, and $a_0$ as before with the additional property that $\norm{a - a_0} < \varepsilon$. Then given $N \in \mathbb{N}$ such that $\norm{\mu^n a} < \varepsilon$ for any $n \geq N$, we have
    \[ \norm{\mu^n a} \leq \norm{\mu^n a_0} + \norm{\mu^n (a - a_0)} < \varepsilon + \varepsilon = 2 \varepsilon. \]
    It follows that $\mu^n a \to 0$. The proof of \thref{statstatesandsimplicity} shows that any $\mu$-stationary state on $C(X) \rtimes_\lambda G$ is of the form we want.
\end{proof}

With \thref{onemeasureconv} in hand, we generalize Hartman and Kalantar's result on C*-simplicity being equivalent to unique stationarity of the action of $G$ on the space of amenable subgroups $\Sub_a(G)$. Recall that Kawabe \cite[Theorem~5.2]{Kawabe} introduced the $G$-space $\Sub_a(X,G)$ of pairs $(x,H)$, where $x \in X$ and $H$ is an amenable subgroup of the stabilizer group $G_x$. 
Observe that the canonical projection $\Sub_a(X,G) \surjectsonto X$ induces an inclusion $C(X) \subseteq C(\Sub_a(X,G))$.

\begin{proof}[Proof of \thref{unique_stationarity_generalized_amenable_subgroup_space}]
    There is a $G$-equivariant, unital and completely positive map $\theta : C(X) \rtimes_\lambda G \to C(\Sub_a(X,G))$ given by $\theta(f \lambda_t)(x,H) = f(x) 1_H(t)$. A similar map can be found used in the proof of \cite[Theorem~5.2]{Kawabe}, but a proof of the existence of such a map is not given. We briefly argue existence here. Given any $(x,H) \in \Sub_a(X,G)$, it is not hard to show that there is a state $\phi \in S(C(X) \rtimes_\lambda G)$ given by $\phi(f \lambda_t) = f(x) 1_H(t)$. This gives us a continuous map from $\Sub_a(X,G)$ to $S(C(X) \rtimes_\lambda G)$, and $\theta : C(X) \rtimes_\lambda G \to C(\Sub_a(X,G))$ is dual to this map.
    
    Choose a generalized measure $\mu \in P(G,C(X))$ as in \thref{statstatesandsimplicity}. Now given a $\mu$-stationary $\eta$ in $P(\Sub_a(X,G))$, we have that $\eta \circ \theta : C(X) \rtimes_\lambda G \to \C$ is necessarily of the form $\nu \circ \bE$. In particular, we note that for $t \neq e$,
    \[ \eta(\{(x,H) | t \in H\}) = \eta(\theta(\lambda_t)) = 0. \]
  Countability of $G$ gives us that $\bigcup_{t \neq e} \{(x,H) | t \in H\}$ is also a null set, or in other words, its complement $X \times \{\{e\}\}$ has measure $1$.
    
    Conversely, assume that the crossed product $C(X) \rtimes_\lambda G$ is not simple. Then by \cite[Theorem~6.1]{Kawabe}, there must exist a closed $G$-invariant subset of $Z \subseteq \Sub_a(X,G)$ that does not intersect $X \times \set{\set{e}}$. Observe that we still canonically have $C(X) \subseteq C(Z)$ by minimality of $X$. Thus, for any $\mu \in P(G,C(X))$, if we choose any $\mu$-stationary state on $C(Z)$ (such a state always exists by \thref{existext}), composing with the canonical quotient $C(\Sub_a(X,G)) \surjectsonto C(Z)$ gives us a $\mu$-stationary state on $C(\Sub_a(X,G))$ with support disjoint from $X \times \set{\set{e}}$.
\end{proof}

\bibliographystyle{amsalpha}
\bibliography{generalized_powers_averaging}

\providecommand{\bysame}{\leavevmode\hbox to3em{\hrulefill}\thinspace}
\providecommand{\MR}{\relax\ifhmode\unskip\space\fi MR }
\providecommand{\MRhref}[2]{%
  \href{http://www.ams.org/mathscinet-getitem?mr=#1}{#2}
}
\providecommand{\href}[2]{#2}
\begin{thebibliography}{BKKO17}

\bibitem[AK20]{AK}
Tattwamasi Amrutam and Mehrdad Kalantar, \emph{On simplicity of intermediate
  {C}*-algebras}, Ergodic Theory and Dynamical Systems \textbf{40} (2020),
  no.~12, 3181--3187.

\bibitem[BK16]{KB}
Rasmus~Sylvester Bryder and Matthew Kennedy, \emph{{Reduced Twisted Crossed
  Products over C*-Simple Groups}}, International Mathematics Research Notices
  \textbf{2018} (2016), no.~6, 1638--1655.

\bibitem[BKKO17]{breuillard_kalantar_kennedy_ozawa_c_simplicity}
Emmanuel Breuillard, Mehrdad Kalantar, Matthew Kennedy, and Narutaka Ozawa,
  \emph{{C*-simplicity and the unique trace property for discrete groups}},
  Publications math{\'e}matiques de l'IH{\'E}S \textbf{126} (2017), no.~1,
  35--71.

\bibitem[Fur63]{Furstenberg}
Harry Furstenberg, \emph{A {P}oisson formula for semi-simple {L}ie groups},
  Annals of Mathematics \textbf{77} (1963), no.~2, 335--386.

\bibitem[Fur73]{Furstenberg1973}
\bysame, \emph{Boundary theory and stochastic processes on homogeneous spaces},
  Harmonic Analysis on Homogeneous Spaces (Calvin~C. Moore, ed.), Proceedings
  of Symposia in Pure Mathematics, vol.~26, American Mathematical Society,
  Providence, Rhode Island, 1973, pp.~193--229.

\bibitem[Haa16]{Haagerup}
Uffe Haagerup, \emph{A new look at {C}*-simplicity and the unique trace
  property of a group}, Operator Algebras and Applications (Cham) (Toke~M.
  Carlsen, Nadia~S. Larsen, Sergey Neshveyev, and Christian Skau, eds.),
  Springer International Publishing, 2016, pp.~167--176.

\bibitem[HK17]{HartKal}
Yair {Hartman} and Mehrdad {Kalantar}, \emph{{Stationary C*-dynamical
  systems}}, arXiv e-prints (2017), arXiv:1712.10133, to appear in Journal of
  the European Mathematical Society.

\bibitem[{Kaw}17]{Kawabe}
Takuya {Kawabe}, \emph{{Uniformly recurrent subgroups and the ideal structure
  of reduced crossed products}}, arXiv e-prints (2017), arXiv:1701.03413.

\bibitem[{Ken}15]{Kennedy}
Matthew {Kennedy}, \emph{{An intrinsic characterization of C*-simplicity}},
  arXiv e-prints (2015), arXiv:1509.01870, to appear in Annales Scientifiques
  de l'{\'E}cole Normale Sup{\'e}rieure.

\bibitem[KK17]{kalantar_kennedy_boundaries}
Mehrdad Kalantar and Matthew Kennedy, \emph{{Boundaries of reduced C*-algebras
  of discrete groups}}, Journal f{\"u}r die reine und angewandte Mathematik
  (Crelles Journal) \textbf{2017} (2017), no.~727, 247--267.

\bibitem[KS19]{KS19}
Matthew Kennedy and Christopher Schafhauser, \emph{Noncommutative boundaries
  and the ideal structure of reduced crossed products}, Duke Mathematical
  Journal \textbf{168} (2019), no.~17, 3215--3260.

\bibitem[Nag20]{Naghavi}
Zahra Naghavi, \emph{{Furstenberg Boundary of Minimal Actions}}, Integral
  Equations and Operator Theory \textbf{92} (2020), no.~2.

\end{thebibliography}
\end{document}